# Entropy and mixing for amenable group actions

By Daniel J. Rudolph and Benjamin Weiss*

## Abstract

For $\Gamma$ a countable amenable group consider those actions of $\Gamma$ as measure-preserving transformations of a standard probability space, written as $\{T_\gamma\}_{\gamma \in \Gamma}$ acting on $(X, \mathcal{F}, \mu)$. We say $\{T_\gamma\}_{\gamma \in \Gamma}$ has *completely positive entropy* (or simply cpe for short) if for any finite and nontrivial partition $P$ of $X$ the entropy $h(T, P)$ is not zero. Our goal is to demonstrate what is well known for actions of $\mathbb{Z}$ and even $\mathbb{Z}^d$, that actions of completely positive entropy have very strong mixing properties. Let $S_i$ be a list of finite subsets of $\Gamma$. We say the $S_i$ *spread* if any particular $\gamma \neq \mathrm{id}$ belongs to at most finitely many of the sets $S_i S_i^{-1}$.

THEOREM 0.1. *For $\{T_\gamma\}_{\gamma \in \Gamma}$ an action of $\Gamma$ of completely positive entropy and $P$ any finite partition, for any sequence of finite sets $S_i \subseteq \Gamma$ which spread we have*
$$\frac{1}{\#S_i} h(\underset{\gamma \in S_i}{\vee} T_{\gamma^{-1}}(P)) \underset{i}{\to} h(P).$$

The proof uses orbit equivalence theory in an essential way and represents the first significant application of these methods to classical entropy and mixing.

## 1. Introduction

The goal of this work is to lift a part of the theory of $K$-actions in the class of measure-preserving transformations of standard probability spaces to actions of countable amenable groups in particular to show that they must be mixing and in fact mixing of all orders. The standard proof of this for transformations (actions of $\mathbb{Z}$ in our vocabulary) or more generally for actions of $\mathbb{Z}^d$ relies on the existence of the *past-algebra* and *tail-algebra* of a partition relative to the action. It seems there is no such good notion for actions of general groups, in particular discrete amenable groups. There are two reasons

*Both authors thank the Institute for Advanced Studies at the Hebrew University of Jerusalem for support of this work.

for working with actions of countable amenable groups. First there is a rather good analogue of the entropy theory of $\mathbb{Z}$ actions for them. This is found in [7] and this as well as the general formulation there of a Rokhlin lemma for such actions are the technical background ideas we will use. Second, these actions are characterized in [1] as precisely those which are orbit-equivalent to actions of $\mathbb{Z}$. We use such an orbit equivalence in an essential way to transfer our problem from $\Gamma$ to $\mathbb{Z}$. Neither entropy nor any mixing properties are orbit equivalence invariants. On the face of it this would seem to block the use of such an orbit equivalence in this area. We avoid this seeming obstacle by working with entropy and mixing properties relative to a sub $\sigma$-algebra and with orbit equivalences that are measurable with respect to this sub $\sigma$-algebra. We show the former are invariants of the latter. This is the first significant application of orbit equivalence theory to classical entropy and mixing questions.

We point out that J. Kieffer [5] constructed a version of the Shannon-McMillan Theorem for amenable group actions predating the work of [7] by the introduction of a generalized tail-field. T. Ward and Q. Zhang in [9] use this method to obtain a conditional entropy theory for such actions. This approach does not appear able to directly give Theorem 2.3 via this generalized tail-field. The work of Ward and Zhang does give essentially all the results we need for conditional entropy. We give an alternative approach to these results anyway building on the methods of [7] rather than those of [5].

In the next section we give a detailed discussion of our results and the proof of Theorem 2.3 except for two results (Theorems 2.6 and 2.11) that require rather extended argument.

In Section 3 we lay out the theory of Rokhlin towers for actions of countable amenable groups, as developed in [7] and the basic geometry of transference leading to Theorem 2.11 on spread sets and a dual result on invariant sets. These results are the tools that allow us to pull entropy and mixing properties through an orbit equivalence. In Section 4 we demonstrate the basic entropy theory we will need to prove Theorem 2.6. As promised we also include the simple argument that the direct product of a cpe action and a Bernoulli action is relatively cpe over the Bernoulli coordinate. This will complete all the technicalities for the proof of Theorem 2.3. In the final section we discuss a number of natural extensions of our work here. These include generalizing our work here to actions of continuous groups and several questions concerning amenable groups that might be answerable using orbit equivalence methods.

## 2. Statements of results and proof of main theorem

Let $\Gamma$ be a countable discrete amenable group. For our purposes the most natural meaning for this is that $\Gamma$ possesses a Følner sequence of sets $H_i$. We will go into more detail on this later. Suffice it here to say that these provide a good analogue for intervals, squares, or boxes of indices in $\mathbb{Z}^d$ for use as both averaging sets for ergodic theorems and for index sets for Rokhlin-like lemmas. Suppose $T = \{T_\gamma\}_{\gamma \in \Gamma}$ is an action of $\Gamma$ by measure-preserving transformations on a standard probability space $(X, \mathcal{F}, \mu)$. If $P$ is a finite partition of $X$ one can define the entropy of the *process* $(T, P)$ as

$$h(T,P) = \lim_{i \to \infty} \frac{1}{\#H_i} h(\bigvee_{\gamma \in H_i} T_{\gamma^{-1}}(P)).$$

That this limit exists and is independent of the choice of Følner sequence is a part of the standard machinery found in [7].

*Definition* 2.1.   We say a $\Gamma$-action $T$ has *completely positive entropy* (or simply is cpe for short) if for any nontrivial finite partition $P$,

$$h(T,P) > 0.$$

This notion is analogous to one of the many definitions of $K$-ness of a single measure-preserving transformation. As said earlier we use this more precise vocabulary to avoid confusion among all these definitions. Our goal is to demonstrate that if $T$ is cpe then it must be mixing of all orders. We will demonstrate something somewhat stronger, that there is a uniformity in this that provides a condition equivalent to cpe. We now formulate this. In a countable group $\Gamma$ we say a sequence of sets $S_i$ *spreads* if each $\gamma \neq \text{id}$ belongs to at most finitely many of the sets $S_i S_i^{-1}$ and for a finite set $K \subseteq \Gamma$ we say $S$ is $K$-spread if for any $\gamma_i \neq \gamma_2 \in S$ we have $\gamma_1 \gamma_2^{-1} \notin K$. If $\Gamma$ is $\mathbb{Z}$, for a sequence of sets to spread simply means the gaps between consecutive terms in the sets grow. The following result is well known for cpe actions of a single transformation and is what we will generalize here.

THEOREM 2.2.   *For $T$ a measure-preserving transformation, $T$ is a $K$-automorphism (is cpe) if and only if for any finite partition $P$ and $\varepsilon > 0$ there is an $N$ so that for any finite subset $S = \{s_1 < s_2 < \ldots < s_t\} \subseteq \mathbb{Z}$ with $s_{i+1} - s_i > N$ for all $i$,*

$$\left|\frac{1}{\#S} h(\bigvee_{\gamma \in S} T_{\gamma^{-1}}(P) - h(P))\right| < \varepsilon.$$

In Theorem 2.2 the conclusion can be written in a couple of equivalent forms. First a form that is *a priori* weaker:

$$\frac{1}{\#S} h(\vee_{j \in S} T^{-j}(P)) > h(P) - \varepsilon$$

and now a form that is *a priori* stronger:

$$\overline{d}_t(\vee_{j \in S} T^{-j}(P); B_t(P)) < \varepsilon$$

where $B_t(P)$ is an independent and identically distributed sequence of $t$ random variables each with the distribution of $P$. It is quite standard to show that these three are all in fact equivalent. The argument is on finite lists of random variables and hence remains true for actions of countable groups. For actions of a single transformation $T$ one can state another equivalent condition, that the sequence of processes $(T^n, P)$ converge in $\overline{d}$ to $B(P)$, the independent and identically distributed process with time-zero distribution that of $P$.

This last formulation may not make sense for an arbitrary discrete amenable group, as it may not have spread co-finite subgroups. But if it does we will be able to reach this conclusion as well.

Here is the goal of our work:

THEOREM 2.3. *For $T$ a cpe action of the discrete amenable group $\Gamma$ and any finite partition $P$ and $\varepsilon > 0$ there is a finite subset $K \subseteq \Gamma$ so that for any finite set $S$ that is $K$-spread*

$$\left| \frac{1}{\#S} h(\vee_{\gamma \in S} T_{\gamma^{-1}}(P)) - h(P) \right| < \varepsilon.$$

To demonstrate Theorem 2.3 we "transfer" the problem to actions of $\mathbb{Z}$ by the fundamental result of [1], that any discrete amenable group action is orbit-equivalent to an action of $\mathbb{Z}$ and by our development here of a relative entropy theory for actions of $\Gamma$.

*Definition* 2.4. For $T$ a measure-preserving action of $\Gamma$, $\mathcal{A}$ a $T$-invariant sub-$\sigma$-algebra and $P$ a finite partition we define the conditional entropy of the process $(T, P)$ given the algebra $\mathcal{A}$ to be

$$h(T, P | \mathcal{A}) = \inf\{h(T, P \vee Q) - h(T, Q) : Q \text{ is } \mathcal{A}\text{-measurable}\}.$$

This is of course just the standard definition for actions of $\mathbb{Z}$. We can also consider orbit equivalences that are $\mathcal{A}$-measurable.

*Definition* 2.5. Let $T$ be a measure-preserving and free action of $\Gamma$ and $\mathcal{A}$ a $T$-invariant sub-$\sigma$-algebra. Suppose $S$ is a free action of some perhaps different group $\Gamma'$ but having the same orbits as $T$. We say *the orbit change from $T$ to $S$ is $\mathcal{A}$-measurable* if for each $\gamma' \in \Gamma'$ the functions $\gamma(x)$ defined by

$$T_{\gamma(x)}(x) = S_{\gamma'}(x)$$

are $\mathcal{A}$-measurable.

For example, suppose $T$ is actually a direct product of two free $\Gamma$-actions $T_1 \times T_2$ and there is an action $U$ of some $\Gamma'$ with the same orbits as $T_2$. The action of $U$ then can be lifted uniquely to an action $\hat{U}$ with the same orbits as $T_1 \times T_2$. If $U_{\gamma'}(x_2) = T_{2,\gamma(x_2)}(x_2)$ then we can (and must) set $\hat{U}_{\gamma'}(x_1, x_2) = (T_1 \times T_2)_{\gamma(x_2)}(x_1, x_2)$. In this case it is clear that the orbit change from $T_1 \times T_2$ to $\hat{U}$ is measurable with respect to the second coordinate $\sigma$-algebra $\mathcal{B}_2$.

Note also that if the orbit change from $T$ to $U$ is $\mathcal{A}$-measurable then $\mathcal{A}$ is a $U$-invariant sub-$\sigma$-algebra.

With these comments in mind we can now state the core technical fact we will use to obtain Theorem 2.3.

THEOREM 2.6. *Suppose $T$ is a free and ergodic action of a countable and discrete amenable group $\Gamma$ and $\mathcal{A}$ is a $T$-invariant sub-$\sigma$-algebra. Suppose also that $U$ is a free (and necessarily ergodic) action of $\Gamma'$ with the same orbits as $T$ ($\Gamma'$ is necessarily amenable). Suppose the orbit change from $T$ to $U$ is $\mathcal{A}$-measurable. Then for any finite partition $P$ we conclude*

$$h(T, P|\mathcal{A}) = h(U, P|\mathcal{A}).$$

Notice that this result is in sharp contrast to the absolute case. All ergodic actions of $\mathbb{Z}$ are orbit-equivalent so entropy itself cannot be invariant. Here that we deal with conditional entropy and require the orbit change to be measurable with respect to the sub-$\sigma$-algebra on which we condition changes the situation entirely. In this situation the orbit change no longer has the freedom to modify the conditional entropy, that is to say the entropy that $T$ has above and beyond the entropy on the sub-$\sigma$-algebra.

With this result in hand and $T$ some cpe action of $\Gamma$ we "transfer" the problem to $\mathbb{Z}$ by considering a direct product action

$$T \times T_2 = \{T_\gamma \times T_{2,\gamma}\}_{\gamma \in \Gamma},$$

where all we really need of $T_2$ is that this direct product should still be ergodic and the action of $T_2$ should be free and hence orbit-equivalent to an action of $\mathbb{Z}$. As described above this means the orbit equivalence of $T_2$ to some $\mathbb{Z}$-action $U$ will lift to an ergodic map $\hat{U}$ giving a $\mathbb{Z}$ action orbit-equivalent to $T \times T_2$. Moreover the orbit change from $T \times T_2$ to $\hat{U}$ is measurable with respect to its second coordinate.

In the direct product of $T$ and $T_2$ consider the relative Pinsker algebra over the second coordinate, that is to say the $\sigma$-algebra generated by all partitions $P$ with $h(T \times T_2, P|\mathcal{B}_2) = 0$. When this algebra is just $\mathcal{B}_2$ the action is said to be relatively cpe over $\mathcal{B}_2$. This relative Pinsker algebra contains both $\mathcal{B}_2$ and the Pinsker algebra of $T$. It is a general fact in a direct product like this that the relative Pinsker algebra is precisely the span of these two algebras. (see [2]). When $T_2$ is a Bernoulli action there is a particularly easy proof so we include this (see Corollary 4.11). In our case this means that the Pinsker algebra of $T \times T_2$ is just that of $T_2$ which is to say $T \times T_2$ is relatively cpe over $\mathcal{B}_2$.

We conclude that $\hat{U}$ must be relatively cpe (relatively $K$) over the algebra $\mathcal{B}_2$ as well. The theory of relative $K$-ness for actions of $\mathbb{Z}$ is well-developed and parallels completely the nonrelative case. We will quite regularly write expressions of the form $E(P|\mathcal{A})$ where $P$ is a finite partition and $\mathcal{A}$ is a sub-$\sigma$-algebra. Consider $P$ as a map from $X$ to the labels $\{1, 2, \ldots, n\}$. By $E(P|\mathcal{A})$ we will mean the probability vector

$$\{E(\mathbf{1}_{P^{-1}(1)}|\mathcal{A}), \ldots, E(\mathbf{1}_{P^{-1}(n)}|\mathcal{A})\}.$$

For $U$ an ergodic measure-preserving transformation, $\mathcal{A}$ a $U$-invariant sub-$\sigma$-algebra, $P$ a finite partition, and $S$ a finite subset of $\mathbb{Z}$ we define

$$L(\vee_{i \in S} U^{-i}(P)|\mathcal{A}) = H(E(\vee_{i \in S} U^{-i}(P)|\mathcal{A})) = E(I(\vee_{i \in S} U^{-i}(P)|\mathcal{A})|\mathcal{A}).$$

That is to say, we calculate the conditional probabilities given $\mathcal{A}$ of the various sets in $\vee_{i \in S} U^{-i}(P)$ and then take the entropy of this probability vector. (For our purposes $h$ calculates the entropy of a partition relative to a fixed measure and $H$ calculates the entropy of a probability vector.)

The following result, due to M. Rahe [8] is fundamental to obtaining our result. Notice that when $\mathcal{A}$ is the trivial algebra it reduces to Theorem 2.2.

THEOREM 2.7 (M. Rahe [8]). *Suppose $\hat{U}$ is an ergodic measure-preserving transformation, $\mathcal{A}$ is a $\hat{U}$-invariant sub-$\sigma$-algebra, and $\hat{U}$ is relatively $K$ with respect to $\mathcal{A}$. Then for any finite partition $P$ and $\varepsilon > 0$ there is an $N$ so that for any finite subset $S \subseteq \mathbb{N}$ with $|s_i - s_j| > N$ for all $i \neq j$,*

$$\left\| \frac{1}{\#S} L(\vee_{i \in S} \hat{U}^{-i} P | \mathcal{A}) - \frac{1}{\#S} \sum_{i \in S} L(\hat{U}^{-i}(P)|\mathcal{A}) \right\|_1 < \varepsilon.$$

Notice that in the conclusion all the information functions are $\mathcal{A}$-measurable but as stated the set $S$ is a constant over the space. We need to consider the possibility that $S = S(x)$ is an $\mathcal{A}$-measurable choice of the set of indices along which the calculation is made. We also will need to loosen up the condition that *all* gaps in the sequence are at least $N$.

*Definition* 2.8. Suppose $S(x)$ is a measurable choice of a $k$-element subset of $\mathbb{Z}$. We say $S(x)$ is $N$-quasi-spread if for all $x$ outside a subset of measure less than $1/N$, there is a subset $S'(x) \subseteq S(x)$ with $\#S'(x)/\#S(x) > 1 - 1/N$ and for all distinct $s, s' \in S'(x)$ we have

$$|s - s'| \geq N.$$

We can generalize this idea to a countable discrete group $\Gamma$ by fixing a listing of its elements $\gamma_1, \gamma_2, \ldots = \Gamma$. Suppose $S(x) = \{s_1(x), \ldots, s_k(x)\}$ is a measurable choice of a $k$-element subset of $\Gamma$. We say $S(x)$ is $N$-quasi-spread if for all $x$ outside a subset of measure less than $1/N$, there is a subset

$S'(x) \subseteq S(x)$ with $\#S'(x)/\#S(x) > 1 - 1/N$ and for all distinct $s, s' \in S'(x)$ we have
$$s^{-1}s' \notin \{\gamma_1, \gamma_2, \ldots, \gamma_N\}.$$

Although this latter definition depends explicitly on how we choose to list the elements of the group it is easy to see that if we fix two listings of the group, to be well quasi-spread for one listing is to be well quasi-spread for the other.

We conjecture that just assuming that $S(x)$ is sufficiently quasi-spread will not be enough to prove Theorem 2.7. Our proof requires us to ask something else of the function that is automatically true for a constant set.

*Definition* 2.9.    For $S(x)$ a measurable choice of $k$-element subset of $\Gamma$ we say $S$ is *uniform* relative to the $\Gamma$-action $T$ if there are $k$ elements $W_i = T_{s_i(x)}(x)$ in the full-group of $T$ (i.e. the $W_i$ are one-to-one, onto and measure-preserving) such that
$$S(x) = \{s_1(x), \ldots, s_k(x)\}.$$

Notice that this implies the $k$ image points $W_i(x)$ are distinct.

The following simple observation is one piece of our transference, saying that the notion of uniformity transfers through a rearranging of the orbit.

LEMMA 2.10.    *Suppose $T$ is a free and ergodic action of the countable group $\Gamma$ and $U$ is an action of another countable group $\Gamma'$ with the same orbits as $T$. Suppose $S(x)$ is a $k$-element set-valued function of $X$ taking values in $\Gamma$ and uniform. Let $W_i(x)$ be $k$ maps whose set of images form the $S(x)$. Each $W_i$ can be written as*
$$W_i(x) = U_{v_i(x)}(x).$$
*Letting $V(x) = \{v_1(x), v_2(x), \ldots, v_k(x)\}$, the set-valued function $V(x)$ takes values in the $k$-point subsets of $\Gamma'$ and is also uniform.*

This piece of the transference is transparent. The other piece will require proof.

THEOREM 2.11.    *Suppose $T$ is a free action of the countable group $\Gamma = \{\gamma_1, \gamma_2, \ldots\}$ and $U$ is a free action of $\Gamma' = \{\gamma'_1, \gamma'_2, \ldots\}$ with the same orbits. Given any $N \in \mathbb{N}$ there is an $M \in \mathbb{N}$ so that for any $k$ and $S(x)$ taking values in the $k$-point subsets of $\Gamma$ that is uniform and $M$-quasi-spread, the set-valued function $V(x)$, taking values in the $k$-point subsets of $\Gamma'$ will be uniform and $N$-quasi-spread.*

We postpone the proof of Theorem 2.11 to Section 2.

In Theorem 2.7 it is enough to ask that the set $S$ be $N$-quasi-spread as the indices outside $S'$ will contribute only an error on the order of $\log \#P/N$ to each term in the difference. What is not at all clear though is that we can replace a constant set $S$ with a variable $S(x)$. We will see though that under the assumption of uniformity we can.

THEOREM 2.12.   *Suppose $U$ is an ergodic measure-preserving transformation, $\mathcal{A}$ is a $U$-invariant sub-$\sigma$-algebra, and $U$ is relatively $K$ with respect to $\mathcal{A}$. Then for any finite partition $P$ and $\varepsilon > 0$ there is an $N$ so that for any $\mathcal{A}$-measurable function $S$ taking values in the $k$-point subsets of $\mathbb{Z}$ that is both $N$-quasi-spread and uniform we will have*

$$\left\| \frac{1}{k} L(\vee_{i \in S} U^{-i} P | \mathcal{A}) - \frac{1}{k} \sum_{i \in S} L(U^{-i}(P)|\mathcal{A}) \right\|_1 < \varepsilon.$$

*Proof.* To say that $U$ is relatively $K$ with respect to the $\sigma$-algebra $\mathcal{A}$ is to say that any partition $P$ has a trivial relative tail. That is to say,

$$\bigcap_N \left( \bigvee_{i \geq N} U^{-i}(P) \vee \mathcal{A} \right) = \mathcal{A}$$

and hence

$$\lim_{N \to \infty} E(P | \vee_{i \geq N} U^{-i}(P) \vee \mathcal{A}) = E(P|\mathcal{A}).$$

Assume $P$ is an $n$-set partition

$$B_N = \{x : |L(P| \vee_{i \geq N} U^{-i}(P) \vee \mathcal{A})(x) - L(P|\mathcal{A})(x)| < \varepsilon/(4 \log n)\}.$$

Choose $N$ sufficiently large that $\mu(B_N) > 1 - \varepsilon/(4 \log n)$ and that $1/N < \varepsilon/(4 \log n)$.

Now assume $S(x)$, taking values in the $k$-point subsets of $\mathbb{Z}$, is both $N$-quasi-spread and uniform. For each $x \in X$ rewrite $S(x)$ as $s_1(x) > s_2(x) > \cdots > s_{t(x)}(x), s_{t(x)+1}(x), \ldots, s_k(x)$ where the ordering is chosen so that $t(x)$ is maximal for the properties

1) for $i < t(x)$, $s_i(x) - s_{i+1}(x) > N$, and

2) for $i \leq t(x)$, $U^{s_i(x)}(x) \in B_N$.

Notice that as $S(x)$ is $\mathcal{A}$-measurable as is $B_N$, we can assume that both $s_i(x)$ and $t(x)$ are as well.

As $S$ is uniform we can calculate

$$\int \#\{s_i(x) : U^{s_i(x)}(x) \in B_N\} \, d\mu \geq k(1 - \varepsilon/(4 \log n))$$

and letting $S'(x) \subseteq S(x)$ be a maximal set of indices pairwise more than $N$ apart, as $S$ is $N$-quasi-spread

$$\int \#S'(x) \, d\mu \geq k(1 - \varepsilon/(2 \log n)).$$

We conclude that
$$\int t(x)\,d\mu \geq k(1 - 3\varepsilon/(4\log n)).$$

We now calculate $L(\vee_{i\in S(x)} U^{-i}(P)|\mathcal{A})$ one term at a time in the form
$$\sum_{i=1}^{k} L(U^{s_i(x)}(P)| \vee_{j<i} U^{s_j(x)}(P) \vee \mathcal{A}).$$

We immediately conclude that
$$L(\vee_{i\in S(x)} U^{-i}(P)|\mathcal{A}) \leq \sum_{i=1}^{k} L(U^{-s_i(x)}(P)|\mathcal{A})$$
and so we only need an inequality in the other direction.

The following inequalities are true almost surely in $x$.
$$t\sum_{i=1}^{k} L(U^{s_i(x)}(P)| \vee_{j<i} U^{-s_j(x)}(P) \vee \mathcal{A})$$
$$\geq \sum_{i=1}^{t(x)} L(U^{s_i(x)}(P)| \vee_{j>s_i(x)+N} U^{-j}(P) \vee \mathcal{A})$$
$$\geq \sum_{i=1}^{t(x)} (L(U^{s_i(x)}(P)|\mathcal{A}) - \varepsilon/(4\log n))$$
$$\geq \sum_{i=1}^{k} (L(U^{s_i(x)}(P)|\mathcal{A}) - \varepsilon/(4\log n)) - (k-t(x))\ln n.$$

We conclude that
$$\int \left| \frac{1}{k} L(\vee_{i\in S(x)} U^{-i}(P)|\mathcal{A}) - \frac{1}{k}\sum_{i\in S} L(U^{-i}(P)|\mathcal{A})\right| d\mu$$
$$\leq \varepsilon/(4\log n) + k - \int t\, d\mu \leq \varepsilon. \qquad \square$$

Combining Theorems 2.12 and 2.11 we reach the following conclusion.

THEOREM 2.13. *Suppose $T$ is a free and ergodic action of a countable amenable group $\Gamma = \{\gamma_1, \gamma_2, \ldots\}$, $\mathcal{A}$ is a $T$-invariant sub-$\sigma$-algebra restricted to which the action of $T$ is still free and relative to which the action of $T$ is cpe, and lastly $P$ is a finite partition. Given any $\varepsilon > 0$ there is an $M$ so that for any $k$ and any $\mathcal{A}$-measurable function $S(x)$ taking values in the $k$-point subsets of $\Gamma$, if $S$ is both uniform and $M$-quasi-spread then we will have*
$$\left\| \frac{1}{k} L(\vee_{\gamma\in S} T_{\gamma^{-1}} P|\mathcal{A}) - \frac{1}{k}\sum_{\gamma\in S} L(T_{\gamma^{-1}}(P)|\mathcal{A}) \right\|_1 < \varepsilon.$$

*Proof.* Construct an action $U$ of $\mathbb{Z}$ that has the same orbits as $T$ and for which the orbit rearrangement is $\mathcal{A}$-measurable. For any $S(x)$ we now obtain a $V(x)$ taking values in the $k$-point subsets of $\mathbb{Z}$ that must be uniform. Notice that

$$L(\vee_{\gamma \in S(x)} T_{\gamma^{-1}}(P)|\mathcal{A})(x) = L(\vee_{i \in V(x)} U^{-i}(P)|\mathcal{A})(x)$$

and for almost every (a.e.) $x$, for each $\gamma \in S(x)$ there is a unique $i \in V(x)$ with $T_\gamma(x) = U^i(x)$ and so

$$L(T_{\gamma^{-1}}(P)|\mathcal{A}) = L(U^{-i}(P)|A).$$

Thus we know

$$\left\| \frac{1}{k} L(\vee_{\gamma \in S} T_{\gamma^{-1}}(P)|\mathcal{A}) - \frac{1}{k} \sum_{\gamma \in S} L(T_{\gamma^{-1}}(P)|\mathcal{A}) \right\|_1$$
$$= \left\| \frac{1}{k} L(\vee_{i \in V} U^{-i}(P)|\mathcal{A}) - \frac{1}{k} \sum_{i \in V} L(U^{-i}(P)|\mathcal{A}) \right\|_1.$$

From Theorem 2.11 we know that for any $N$ if $M$ is large enough then $V$ will be $N$-quasi-spread and Theorem 2.12 completes the theorem. □

In the case we are interested in $T$ is in fact $T_1 \times T_2$, the algebra $\mathcal{A}$ is the second coordinate algebra, $P$ is measurable with respect to the first coordinate algebra and $S$ is a constant choice of set. In this case notice that

$$L(\vee_{\gamma \in S} T_{\gamma^{-1}}(P)|\mathcal{A}) = h(\vee_{\gamma \in S} T_{\gamma^{-1}}(P))$$

and

$$L(T_{\gamma^{-1}}(P)|\mathcal{A}) = h(P).$$

From this the proof of Theorem 2.3 is complete.

## 3. Discrete amenable groups and spread sets

In this section we lay out the basic ergodic theory of measure-preserving actions of countable and discrete amenable groups. This material can all be found in [7] or [4]. Our goal is to describe some particular structures that are invariants of orbit equivalence. As described in the introduction we will be "transferring" our work from an action of $\Gamma$ to one of $\mathbb{Z}$ with the same orbits. We here establish a vocabulary for transferring certain calculations through this transference. The most important of these will be conditional entropy. The transference of these calculations will be made by showing that they can be made on towers. Our first step is to lay out the structure of Rokhlin towers as developed in [7]. These build from notions of invariance of sets in $\Gamma$.

Notions of essential invariance of finite subsets $F \subseteq \Gamma$ is central to the study of the ergodic theory of discrete amenable group actions. Here is the fundamental definition used in [7].

*Definition* 3.1. Let $\delta > 0$. Let $K \subseteq \Gamma$ be a finite set. A subset $F \subseteq \Gamma$ is called $(\delta, K)$-*invariant* if
$$\frac{\#(KK^{-1}F \triangle F)}{\#F} < \delta.$$

To say that a set $F$ is *sufficiently invariant* means that for some (unspecified) $\delta > 0$ and finite $K \subseteq \Gamma$, $F$ is $(\delta, K)$-invariant.

To say that a list of sets $F_1, F_2, \ldots, F_k$ is *sufficiently invariant* is to say that for some $\delta > 0$ and finite set $K \subseteq \Gamma$, setting $F_0 = K$, for each $j \in \{1, 2, \ldots, k\}$, the set $F_j$ is $(\delta, F_{j-1})$-invariant.

The existence of $(\delta, K)$-invariant sets for all $\delta > 0$ and finite sets $K$ is equivalent to the amenability of the countable group $\Gamma$.

For purely technical reasons a different but equivalent calculation of the degree of invariance of a set is more convenient for us.

*Definition* 3.2. Let $\delta > 0$ and $K \subseteq \Gamma$ be a finite set. A subset $F \subseteq \Gamma$ is called $[K, \delta]$-*invariant* if
$$\#\{\gamma \in F : KK^{-1}\gamma \subseteq F\} > (1 - \delta)\#F.$$

One calculates that if $F$ is $(\delta/(\#K)^2, K)$-invariant, then $F$ is $[K, \delta]$-invariant and conversely if $F$ is $[K, \delta/(\#K)^2]$-invariant, then $F$ is $(\delta, K)$-invariant. Hence when one says that either $F$ or a list $F_1, \cdots, F_k$ are *sufficiently invariant* one need not distinguish which of the two notions, Definition 3.1 or Definition 3.2 is meant.

We now present the Ornstein-Weiss quasi-tiling theorem.

*Definition* 3.3. A finite list of sets $H_1, H_2, \ldots, H_k \subseteq \Gamma$, with $\text{id} \in H_i$, for all $i$, is said to $\varepsilon$-*quasi-tile* a finite set $F \subseteq \Gamma$ if there exist "centers" $c_{i,j}$, $i = 1, 2, \ldots, k$, $j = 1, 2, \ldots, l(i)$, and subsets $H_{i,j} \subseteq H_i$ such that:

1. $\#H_{i,j} \geq (1-\varepsilon)\#H_i$, for $j = 1, \ldots, l(i)$,
2. the $H_{i,j}c_{i,j} \subseteq F$ are disjoint, and
3. $\#(\underset{i,j}{\cup}H_{i,j}c_{i,j}) \geq (1-\varepsilon)|F|$.

THEOREM 3.4 ([7]). *Given $\varepsilon > 0$, there exists $N = N(\varepsilon)$ such that in any countable discrete amenable group $\Gamma$, if $H_1, \ldots, H_N$ is any sufficiently invariant list of sets, then for any $D \subseteq \Gamma$ that is sufficiently invariant (depending on the choice of $H_1, \ldots, H_N$), $D$ can be $\varepsilon$-quasi-tiled by $H_1, \ldots, H_N$.*

This theorem is the essential content of Theorem 6, I.2 [7]. Our definition of $\varepsilon$-quasi-tiling is slightly different; weaker in that we do not ask that $H_i c_{i,j} \cap H_k c_{k,l} = \emptyset, i \neq k$, and stronger in that we require $H_{i,j} c_{i,j} \subseteq F$. Obtaining the latter from Theorem 6, I.2 [7] is easy if $D$ is sufficiently invariant and N is fixed.

Suppose $(X, \mathcal{B}, \mu)$ is a standard probability space. Suppose $T$ is a measure-preserving free action of $\Gamma$ on $X$. It is useful to keep in mind that $T$ is a function of two variables, i.e., is a map from $\Gamma \times X$ to $X$. For a finite set $F \subseteq \Gamma$ and measurable subset $A \in \mathcal{B}$ with $\mu(A) > 0$, consider $F \times A \subseteq \Gamma \times X$. As a measure on $F \times A$, put the direct product $c \times \mu$ of counting measure $c$ and $\mu$. Consider the map $T$ restricted to this rectangle $F \times A$. On each level set $g \times A$, $T$ is one-to-one and measure-preserving. On any fiber set $F \times x$, $T$ is again one-to-one. We definitely do not expect $T$ to be one-to-one on $F \times A$. It is clear, though, that $T$ is nonsingular and, at most, $\#F$ to 1.

For $S \subseteq F \times A$, let

$$S(x) = \#\{g \in F; (g, x) \in S\}$$

and

$$c_S(x) = \#S(x).$$

Of course

$$(c \times \mu)(S) = \int_A c_S(x) \, d\mu(x).$$

Set

$$c(S) = \min_{x \in A} c_S(x).$$

*Definition* 3.5. We say that $F \times A$ *maps an $\varepsilon$-quasi-tower* if there exists a measurable subset $S \subseteq F \times A$ such that:

1. $T|_S$ is one-to-one, and

2. $c(S) \geq (1 - \varepsilon)\#F$.

We say $F \times A$ *maps a real-tower* if $T|_{F \times A}$ is one-to-one.

The $\varepsilon$-quasi-tower itself (or real-tower as the case may be) is the set $T(F \times A) \subseteq X$. Notice that we may always assume $T(S) = T(F \times A)$. Also notice that if there exists an $S \subseteq F \times A$, such that $T$ is one-to-one on $S$ and $(c \times \mu)(S) > (1 - \varepsilon^2)(c \times \mu)(F \times A)$, then there must exist an $A' \subseteq A$, with $\mu(A') > (1 - \varepsilon)\mu(A)$, such that $c_S(x) > (1 - \varepsilon)\#F$, for all $x \in A'$. Hence $F \times A'$ maps to an $\varepsilon$-quasi-tower. We will return to real-towers in a bit but for now we focus on $\varepsilon$-quasi-towers, which are the core of the Rokhlin lemma of [7]. This is only a minor modification of Theorem 5, II.2 of [7].

THEOREM 3.6 ([7]). *Suppose $\Gamma$ is a discrete amenable group. For any $\varepsilon > 0$, there exist $\delta > 0$, $K \subseteq \Gamma$ and $N = N(\varepsilon)$ such that for any sequence $H_1, \ldots, H_N$ of $[K, \delta]$-invariant subsets of $\Gamma$, and any free measure-preserving $\Gamma$-action $T = \{T_g\}_{\gamma \in \Gamma}$, acting on $(X, \mathcal{B}, \mu)$, there exist sets $A_1, \ldots, A_N \in \mathcal{B}$ such that*:

1. *each $H_i \times A_i$ maps to an $\varepsilon$-quasi-tower $\mathcal{R}_i$ in $X$,*

2. *for $i \neq j$, $\mathcal{R}_i \cap \mathcal{R}_j = \emptyset$, and*

3. $\mu(\cup_{i=1}^N \mathcal{R}_i) > 1 - \varepsilon$.

A collection of sets of the form $\{H_i \times A_i\}_{i=1}^N$ satisfying (1), and (2) we say forms $\varepsilon$-*Rokhlin towers*. We indicate (3) by saying the towers cover $(1 - \varepsilon)$ of $X$. If the $H_i \times A_i$ are real-towers satisfying (2) then we say $\{H_i, A_i\}$ is a *castle* (again covering $(1 - \varepsilon)$ of $X$.) We will be using castles a lot as we proceed so we stress what they are. A *castle* is a collection of the form $\{H_i, A_i\}$ where each pair $H_i, A_i$ form a real-Rokhlin tower and the list of towers have disjoint images. We refer to the union of these images in $X$ as the *castle image*.

Notice that if $\{H_i \times A_i\}$ form an $\varepsilon$-quasi-Rokhlin tower then in each of the $H_i \times A_i$ there is a subset $S_i$ with $c(S_i) > (1 - \varepsilon)\#H_i$. For each $x \in A$ remember $S_i(x) = \{\gamma : (\gamma, x) \in S_i\}$. Partition $A_i$ into subsets $A_{i,j}$ on which the choice of set $S_i(x)$ is a constant $H_{i,j}$. Now $\#H_{i,j} \geq \#H_i(1-\varepsilon)$. The collection of sets $\{H_{i,j} \times A_{i,j}\}$ forms a castle covering the same set as the original $\varepsilon$-quasi-Rokhlin tower. If each $H_i$ is $(K, \delta)$-invariant, then the new $H_{i,j}$ are still $(K, \delta + 2\varepsilon)$-invariant. Furthermore if the sets $A_i$ were known to be measurable with respect to some invariant sub-$\sigma$-algebra then the bases of the castle could also be selected to be measurable with respect to that invariant algebra.

We now state and prove a castle based version of the mean ergodic theorem for actions of countable amenable groups. This fact is not difficult to prove. We choose a version of the proof which is parallel to our proof of Theorem 2.6 as a warm-up to that proof.

For $T$ an action of $\Gamma$, $f$ a real-valued function, and $\{A_i, H_i\}$ forming towers (real or quasi does not matter yet) we can ask how close averages of $f$ across the towers approximate the integral of $f$, that is to say we can ask how small is

$$\sum_i \int_{A_i} \left| \sum_{\gamma \in H_i} f(T_\gamma(x)) - \#H_i \int f \, d\mu \right| d\mu.$$

It may appear that a normalization by $\#H_i$ is missing but the sum is normalizing each tower to have mass $\#H_i \mu(A_i)$.

The next lemma gives the pivotal observation, that if there exist towers where we are seeing convergence on average, then on all sufficiently invariant towers we must also be seeing such convergence.

LEMMA 3.7. *Suppose $T$ is a free and ergodic action of the countable group $\Gamma$ and $f \in L_1(\mu)$. Given any $\varepsilon > 0$ there is a $\delta > 0$ so that if there exists a castle $\{A_i, H_i\}$ covering all but $\delta$ of $X$ and satisfying*

$$\sum_i \int_{A_i} \left| \sum_{\gamma \in H_i} f(T_\gamma(x)) - \#H_i \int f \, d\mu \right| d\mu < \delta,$$

*then for any sufficiently invariant $K_j$, if a castle $\{B_j, K_j\}$ covers at least $\varepsilon$ of $X$ we will have*

$$\sum_j \int_{B_j} \left| \sum_{\gamma \in K_j} f(T_\gamma(x)) - \#K_j \int f \, d\mu \right| d\mu < \varepsilon.$$

*Proof.* Without loss of generality we can assume that $\int f \, d\mu = 0$ and hence that

$$\sum_i \int_{A_i} \left| \sum_{\gamma \in H_i} f(T_\gamma(x)) \right| d\mu < \delta.$$

Also without loss of generality we can assume that there are only finitely many sets $H_i$ and hence there is some finite $K$ containing all of them. Choose $\delta_1$ so small that for any set $C$ with $\mu(C) < \delta_1$ we must have $\int_C |f| \, d\mu < \varepsilon/2$. Assume all the $K_j$ are $[K, \delta_1/3]$-invariant and be sure $\varepsilon(\delta_1/3)^2 > \delta$. Let $\mathcal{T} \subseteq X$ be the tower image of $\{A_i, H_i\}$ and suppose the castle $\{B_j, K_j\}$ covers at least $\varepsilon$ of $X$.

Consider the set $D_j \subseteq B_j$ where

$$\#\{\gamma \in K_j : T_\gamma(x) \notin \mathcal{T}\} \geq \#K_j \delta_1/3.$$

As $1 - \delta < \mu(\mathcal{T}) \leq 1 - \sum_j \mu(D_j) \#K_j \delta_1/3$, we must have

$$\varepsilon \delta_1/3 \geq \sum_j \mu(D_j) \#K_j,$$

which is to say

$$\sum_j \mu(D_j^c) \#K_j \geq \sum_j \mu(B_j) \#K_j (1 - \delta_1/3) \geq \varepsilon(1 - \delta_1/3).$$

For each $B_j$ and each $x \in D_j^c$ consider the collection of group elements

$$d_j(x) = \{\gamma \in K_j : \text{with } T_\gamma(x) \notin \mathcal{T} \text{ or } K\gamma \text{ is not contained in } K_j\}.$$

For $x \in D_j^c$ we know $\#d_j(x) \leq 2\delta_1 \#K_j/3$.

Define a global "bad" set $D$ in the castle image by

$$D = \bigcup_j \left( \bigcup_{x \in D_j} \{T_\gamma(x) : \gamma \in K_j\} \cup \bigcup_{x \notin D_j} \{T_\gamma(x) : \gamma \in d_j(x)\} \right).$$

We calculate
$$\mu(D) \leq \sum_j \mu(D_j)\#K_j + \mu(B_j)\#K_j 2\delta_1/3 \leq \delta_1.$$

Split $f$ into two functions as $f = f_1 + f_2$ where
$$f_1(x) = \begin{cases} 0 & \text{if } x \in D \\ f(x) & \text{otherwise.} \end{cases}$$

We know
$$\sum_j \int_{B_j} \left| \sum_{\gamma \in K_j} f_2(x) \right| d\mu < \varepsilon/2.$$

To show the same for $f_1$ notice that
$$\begin{aligned}
\sum_j \int_{B_j} \left| \sum_{\gamma \in K_j} f_1(T_\gamma(x)) \right| d\mu &= \sum_j \int_{D_j^c} \left| \sum_{\gamma \in d_j(x)} f(T_\gamma(x)) \right| d\mu \\
&= \sum_{i,j} \int_{D_j^c} \left| \sum_{\substack{\gamma' \in d_j(x) \\ T_{\gamma'}(x) \in A_i}} \left( \sum_{\gamma \in H_i} f(T_{\gamma\gamma'}(x)) \right) \right| d\mu \\
&\leq \sum_{i,j} \int_{D_j^c} \sum_{\substack{\gamma' \in d_j(x) \\ T_{\gamma'}(x) \in A_i}} \left| \left( \sum_{\gamma \in H_i} f(T_{\gamma\gamma'}(x)) \right) \right| d\mu \\
&\leq \sum_i \int_{A_i} \left| \sum_{\gamma \in H_i} f(T_\gamma(x)) \right| d\mu
\end{aligned}$$

as
$$\bigcup_{x \in D_j^c} \left( \bigcup_{\substack{\gamma' \in d_j(x) \\ T_{\gamma'}(x) \in A_i}} T_{\gamma'}(x) \right) \subseteq A_i$$

is a disjoint union. This latter is less than $\varepsilon/2$, finishing the result. □

THEOREM 3.8. *For $T$ an ergodic action of the countable amenable group $\Gamma$ and $f \in L_1(\mu)$ and $\varepsilon > 0$ for any sufficiently invariant castle $\{A_i, H_i\}$ whose image covers at least $\varepsilon$ of $X$ we will have*
$$\sum_i \int_{A_i} \left| \sum_{\gamma \in H_i} f(T_\gamma(x)) - \#H_i \int f \, d\mu \right| d\mu < \varepsilon.$$

*Proof.* From the previous lemma all we need to know is that for any $\delta > 0$ there exist towers covering all but $\delta$ of $X$ and for which the estimate holds with $\delta$ replacing $\varepsilon$. Here is a very simple way to get them. Construct a direct product of $T$ with some auxiliary free and ergodic action $T_1$ of $\Gamma$ so that the direct product $T \times T_1$ is still ergodic. In this product consider only towers $\{B_j, K_j\}$ that are measurable with respect to the second coordinate algebra. As the $\sigma$-algebra generated by the action of $T$ on $f$ is independent of the second coordinate algebra, if we know all the $K_j$ are so invariant that from the mean ergodic theorem

$$\int \left| \sum_{\gamma \in K_j} f(T_\gamma(x)) - \#K_j \int f \, d\mu \right| d\mu < \delta \#K_j$$

then the tower estimate will automatically follow.  □

We will reuse this little trick of adding on a second coordinate in an even more decisive way later in proving Theorem 2.6. The next definition is dual to that of a set being well quasi-spread. For consistency we make the definition more parallel by fixing a listing of the elements of $\Gamma = \{\gamma_1, \gamma_2, \ldots\}$.

*Definition* 3.9.  We say that the castle $\{H_i \times A_i\}$ is *N-invariant* if setting $K = \{\gamma_1, \ldots, \gamma_N\}$

$$\sum_i \#\{\gamma \in H_i : KK^{-1}\gamma \nsubseteq H_i\}\mu(A_i) < 1/N.$$

Notice that if each of the sets $H_i$ is individually $[K, 1/N]$ invariant then the castle will be $N$-invariant, but the invariance of the castle is an average notion and so $1/N^2$ invariance of the castle will imply $[\{\gamma_1, \ldots, \gamma_N\}, 1/N]$-invariance of all but $1/N$ of the $\{A_i, H_i\}$ (measured by their fraction of the castle image).

We have already seen how a spread set transfers through an orbit equivalence to another action. We describe the parallel transference for castles. Suppose $T$ is a free and ergodic action of a countable group $\Gamma$ and $U$ is a free action of $\Gamma'$ with the same orbits. What lies behind the transference is of course the cocycle of the relation $\alpha : X \times \Gamma \to \Gamma'$ where $T_\gamma(x) = U_{\alpha(x,\gamma)}(x)$. It is useful to also give a name to the cocycle in the other direction $\beta : X \times \Gamma' \to \Gamma$ where $U_{\gamma'}(x) = T_{\beta(x,\gamma')}(x)$. Suppose $\{A_i, H_i\}$ is a castle for the action of $T$. For each $A_i$ and point $x \in A_i$ consider the set of points $\mathcal{T}(x) = \{T_\gamma(x) : \gamma \in H_i\}$. This is the slice through the tower image passing through $x$. Setting $K_i(x) = \alpha(x, H_i)$ we have

$$\mathcal{T}(x) = \{U_{\gamma'}(x) : \gamma' \in K_i(x)\}.$$

Partition $A_i$ into subsets $B_{i,j}$ according to the set $K_i(x) = K_{i,j}$. Now the full collection $\{K_{i,j}, B_{i,j}\}$ is a castle for the action of $U$. We refer to this as the image castle. It will be convenient to drop the double-indexing and write these towers simply as $\{K_j, B_j\}$. Notice that if both the orbit change from $T$ to $U$ and the sets $A_i$ are measurable with respect to some invariant sub-$\sigma$-algebra $\mathcal{A}$ then the sets $B_j$ will automatically be $\mathcal{A}$-measurable.

We now verify Theorem 2.11 and the corresponding result for castles. The following result is not essential for our work but we include it as it is perhaps the most primitive of this type of result.

LEMMA 3.10. *Suppose $T$ and $U$ are orbit-equivalent free actions of amenable groups $\Gamma$ and $\Gamma'$ with $\alpha$ and $\beta$ the cocycles as described above. Given any $K' \subseteq \Gamma'$ and $\delta' > 0$ there is a $K \subseteq \Gamma$ and $\delta > 0$ so that if $F \subseteq \Gamma$ is $[K, \delta]$-invariant, then*

$$\mu(\{x : \alpha(x, F) \text{ is } [K', \delta']\text{-invariant}\}) > 1 - \delta'.$$

*Proof.* For $x \in X$ set $\kappa(x) = \beta(x, K'K'^{-1}) \subseteq \Gamma$. Choose the set $K_1$ so that

$$\mu(\{x : \kappa(x) \subseteq K_1\}) > 1 - \delta'/3$$

and set $B = \{x : \kappa(x) \subseteq K_1\}$. Choose $K$ containing $K_1$ and $\delta \leq \delta'/3$ so that if $F$ is $[K, \delta]$-invariant then for all but a set of $x$ of measure at most $\delta'$, the density in the set of points $\{T_\gamma(x) : \gamma \in F\}$ of occurrences of the set $B$ is at least $1 - \delta'/3$ (using the mean ergodic theorem.) Notice that $F$ must be $(K_1, \delta/3)$-invariant. For such an $x$ set $F' = \alpha(x, F)$ and we get a lower bound for $\#\{\gamma'; K'K'^{-1}\gamma' \subseteq F'\}$ by noting it must contain all $\alpha(x, \gamma)$ where $\gamma \in F$, $T_\gamma(x) \in B$ and $K_1\gamma \subseteq F$. Both of these subsets of $F$ omit a fraction of at most $\delta'/3$ from $F$. □

THEOREM 3.11. *Suppose $T$ and $U$ are orbit-equivalent free actions of amenable groups $\Gamma$ and $\Gamma'$ with the cocycles $\alpha$ and $\beta$ as described above. Given any $1 > a > 0$ and $N$, there is an $M$ so that for all castles $\{A_i, H_i\}$ for the action of $T$ that are $M$-invariant and cover $a$ of $X$, the image castle $\{K_j, B_j\}$ relative to the action $U$ must be $N$-invariant.*

*Proof.* We assume $M$ is divisible by $N$ and write it as $M = M'N$. Assuming $\{A_i, H_i\}$ to be $M$-invariant there must be a subset of indices $I$ so that

$$\sum_{i \in I} \mu(A_i) \#H_i > 1 - \frac{1}{2N}\mu \text{ (the castle image)}$$

and for $i \in I$ the set $H_i$ is $(\{\gamma_1, \ldots, \gamma_{M'}\}, 1/M')$-invariant. For $K' = \{\gamma'_1, \ldots, \gamma'_N\}$ choose $K$ so that setting

$$B = \{x : \beta(x, K'K'^{-1})\} \subseteq K$$

we have $\mu(B) > 1 - a(4N^3)^{-2}$. Let $A'_i \subseteq A_i$ consist of those $x$ for which

$$\#\{\gamma \in H_i : T_\gamma(x) \in B\} > (1 - \frac{1}{4N^3})\#H_i.$$

As $\sum_i \mu(A_i) \# H_i > a$ we must have

$$\sum_i \mu(A'_i) \# H_i > a(1 - \frac{1}{4N^3}).$$

Choose $M' > 2N^3$ and so large that $K \subseteq \{\gamma_1, \ldots, \gamma_{M'}\}$. Thus for $i \in I$ $H_i$ is $[K, (2N^2)^{-1}]$-invariant and so for all $i \in I$

$$\#\{\gamma : K\gamma \subseteq H_i\} > (1 - \frac{1}{2N^3}) \# H_i.$$

For $x \in A'_i$ consider $\kappa(x) = \alpha(x, H_i)$. Notice that

$$\begin{aligned}
\#\{\gamma' \in \kappa(x) &: K'K'^{-1} \not\subseteq \kappa(x)\} \\
&= \#\{\gamma = \beta(x, \gamma') : \gamma \in H_i, \beta(x, K'K'^{-1}\gamma') \\
&\quad = \beta(U_{\gamma'}(x), K'K'^{-1})\beta(x, \gamma') \not\subseteq H_i\} \\
&\leq \#\{\gamma' : \beta(U_{\gamma'}(x), K'K'^{-1}) \not\subseteq K\} + \#\{\gamma \in H_i : K\gamma \not\subseteq H_i\} \\
&\leq \#\{\gamma \in H_i : T_\gamma(x) \notin B\} + \frac{\#H_i}{2N^3} \\
&\leq \frac{\#H_i}{N^3}.
\end{aligned}$$

This says for such an $x$ the set $\kappa(x)$ is $[K', 1/N]$-invariant completing the proof. □

The proof of Theorem 2.11 is completely parallel to this and we now present it.

*Proof of Theorem* 2.11. Set $K' = \{\gamma'_1, \ldots, \gamma'_N\}$ and as in the previous proofs choose $K$ so large that

$$B = \{x : \beta(x, K') \subseteq K\}$$

satisfies

$$\mu(B) > 1 - \frac{1}{(4N)^2}.$$

Choose $M$ large enough that $M > 4N$ and so that $K \subseteq \{\gamma_1, \ldots, \gamma_M\}$. Suppose $S : X \to \{$the $k$ point subsets of $\Gamma\}$ is both uniform and $M$-quasi-spread. Write

$$T_{S(x)}(x) = \{V_1(x), \ldots, V_k(x)\}$$

where the $V_i$ are in the full-group of $T$ (and hence are one-to-one and measure-preserving) and have distinct images.

As $\mu(B^c) < (4N)^{-2}$,

$$\mu(\{x : \#\{i \leq k : V_i(x) \notin B\} > \frac{1}{4N}k\}) < \frac{1}{4N}.$$

Set $A = \{x : \#\{i \leq k : V_i(x) \notin B\} \leq k/(4N)\}$. Consider those $x \in A$ for which there is a subset $S'(x) \subseteq S(x)$ occupying a fraction $1 - 1/M$ of $S(x)$ that is $M$-spread. This set has measure at least $1 - 1/M - 1/(4N) > 1 - 1/N$.

For such an $x$ set
$$S''(x) = S'(x) \cap \{\gamma \in S : T_\gamma(x) \in B\}$$
and we note
$$\#S''(x) > \#S'(x) - \frac{k}{4N} > k(1 - \frac{1}{N}).$$

In $V(x) = \alpha(x, S(x))$ set $V'(x) = \alpha(x, S''(x))$. For $\gamma' \in V'(x)$ and $\gamma'_1 \in K'$ with $\gamma'_1 \neq \mathrm{id}$, we need only show that $\gamma'_1 \gamma' \notin V'(x)$ to be finished. This is the same as saying
$$\beta(x, \gamma'_1 \gamma') \notin S''(x).$$

Now
$$\beta(x, \gamma'_1 \gamma') = \beta(U_{\gamma'}(x), \gamma'_1)\beta(x, \gamma').$$

As $\gamma' \in V'(x)$, $U_{\gamma'}(x) \in B$ and $\beta(x, \gamma') \in S'(x)$. As $\gamma' \in K'$, $\beta(U_{\gamma'}(x), \gamma'_1) \in K$. We must conclude, as $S'(x)$ is $M$-spread, that $\beta(U_{\gamma'}(x), \gamma'_1)\beta(x, \gamma')$ cannot be in $S'(x)$ which contains $S''(x)$. □

## 4. Conditional entropy theory and proof of Theorem 2.6

The Shannon-McMillan theorem of [7] states that if $T$ is an ergodic action of $\Gamma$ on $(X, \mathcal{F}, \mu)$ and $P$ is a finite partition then for any $\delta > 0$ if $F \subseteq \Gamma$ is sufficiently invariant then among the atoms of $\bigvee_{\gamma \in F} T_{\gamma^{-1}}(P)$ is a collection $G$ with:

1) $\mu(\cup_{\eta \in G} \eta) > 1 - \delta$,

2) $\#(G) < 2^{(h(T,P)+\delta)\#F}$, and

3) for any $\eta \in G$, $\mu(\eta) < 2^{-(h(T,P)-\delta)\#F}$.

We need a version of this fact relative to, (or conditioned on) a $T$-invariant sub-$\sigma$-algebra. To obtain this requires that we develop some basic conditional entropy theory for actions of countable amenable groups. T. Ward and Q. Zhang have given a development of conditional entropy for amenable group actions but their results (Theorem 3.2 and Corollary 3.3 of [9]) are not quite enough for our needs. Although it is possible to build from their development and get what we need for the sake of consistency and completeness we give the whole story here.

First two primitive observations about entropy that hold equally well for us here. If $T$ is an action of $\Gamma$ and $F$ is some finite subset of $\Gamma$ with

$$Q \subseteq \bigvee_{\gamma \in F} T_{\gamma^{-1}}(P)$$

a finite partition then

$$h(T, Q) \leq h(T, P).$$

Next, for any vector $\vec{p} = \{p_1, \ldots, p_t\}$ of values $p_i \geq 0$ with $\sum p_i \leq 1$ (yes, $\leq$) we can set $p_{t+1} = 1 - \sum p_i$ and let

$$H(\vec{p}) = -\sum_{i=1}^{t+1} p_i \log(p_i).$$

Suppose $P$ and $P'$ are partitions of $X$ labeled by the same symbols, i.e. $P, P' : X \to \Sigma$ where $\Sigma$ is a finite set. We write

$$P \triangle P' = \{x : P(x) \neq P'(x)\}.$$

The following inequality always holds:

$$|h(T, P) - h(T, P')| \leq H(\mu(P \triangle P')) + \mu(P \triangle P') \log(\#\Sigma).$$

In particular among $\Sigma$-valued partitions $h(T, P)$ is uniformly continuous in the $\mu(\cdot \triangle \cdot)$ metric.

*Definition* 4.1.  We abstract the notion in the conclusion of the Shannon-McMillan theorem of [7] as follows. For $K \subseteq \Gamma$ a finite set and $p$ a probability measure on $\Sigma^K$, we say $p$ is $(h, \delta)$-flat if there is a subset $S_0 \subseteq \Sigma^K$ with:

1) $p(S_0) > 1 - \delta$,

2) $\#S_0 \leq 2^{(h+\delta)\#K}$, and

3) for any $\eta \in S_0$, $\mu(\eta) \leq 2^{-(h-\delta)\#K}$.

For a probability measure on $\Sigma^K$ like $p$ in this definition let

$$h_K(p) = -\frac{1}{\#K} \sum_{\eta \in \Sigma^K} p(\eta) \log(p(\eta)).$$

The following are basic calculations.

LEMMA 4.2.  *If $p$ satisfies* 1) *and* 2) *above then*

$$h_K(p) \leq H(\delta) + \delta \log(\#\Sigma) + (1-\delta)(h+\delta) \leq h + H(\delta) + \delta \log(\#\Sigma).$$

*If $p$ satisfies* 1) *and* 3) *above then*

$$h_K(p) \geq (1-\delta)(h-\delta) \geq h - \delta(\log(\#\Sigma) + 1).$$

AMENABLE CPE-ACTIONS MIX1139

Now for the technical fact that if one writes a "flat" distribution on names as a convex combination of distributions without dropping the average entropy by much, then most of the terms in the convex combination must also be pretty flat.

THEOREM 4.3.  *Given $h \geq 0$ and $\delta > 0$ there is a $\delta_1$ so that if $\#K > 1/\delta_1$, $p$ is $(h, \delta_1)$-flat and we write $p$ as a convex combination of measures*

$$p = \int p_t \, dt$$

*with*

$$\int h_K(p_t) \, dt \geq h_K(p) - \delta_1,$$

*then*

$$\mu(\{t : p_t \text{ is } (h, \delta)\text{-flat }\}) > 1 - \delta.$$

*Proof.* In this calculation we will use $F_i(z)$ to represent a sequence of functions, all of which tend to 0 as $z \to 0$. We will find $F_4$ so that for all but $F_4(\delta_1)$ of the values $t$, $p_t$ is $(h, F_4(\delta_1))$-flat.

Assume $K$ and $p$ are as described for $h$ and $\delta_1$. Let $S_0 \subseteq \Sigma^K$ be the set of names making $p$ an $(h, \delta_1)$-flat distribution. For all but $\sqrt{\delta_1}$ of the values $t$,

$$p_t(S_0) \geq 1 - \sqrt{\delta_1},$$

for example, $p_t$ satisfies 1) and 2) of $(h, \sqrt{\delta_1})$-flatness and hence

$$h_K(p_t) \leq h + H(\sqrt{\delta_1}) + \sqrt{\delta_1}(\log(\#\Sigma) + 1).$$

Thus there is a function $F_1$ and for all but $F_1(\delta_1)$ of the values $t$,

$$h_K(p_t) \leq h + F_1(\delta_1).$$

For all $t$ we have $h_K(p_t) \leq \log(\#\Sigma)$ so as $\int h_K(p_t) \, dt \geq h - \delta_1$ there is an $F_2$ and for all but $F_2(\delta_1)$ of the values $t$,

$$h - F_2(\delta_1) \leq h_K(p_t) \leq h + F_2(\delta_1)$$

and

$$p_t(S_0) \geq 1 - F_2(\delta_1).$$

Suppose for such a $p_t$ there is a $c \geq 0$ and a set $S_1 \subseteq \Sigma^K$ with $p_t(S_1) = c$ and for any name $\eta \in C$

$$p_t(\eta) > 2^{-(h-c)\#K}.$$

We want to see that $c$ must be small.

Split $\Sigma^K$ into three sets:

1) $S_1$ of measure $c$,

2) $S_0 \backslash S_1$ of measure $> 1 - \delta_1 - c$, and

3) the rest, of measure $< \delta_1$.

We calculate, conditioning on this three set-partition, that

$$h_K(p_t) \leq \frac{H(c, \delta_1)}{\#K} + c(h-c) + (1 - \delta_1 - c)(h + \delta_1) + \delta_1 \log(\#\Sigma)$$

and

$$h - F_2(\delta_1) \leq h - c^2 + \delta_1 \log(\#\Sigma) + \frac{\log 3}{\#K} + \delta_1$$

and as $\#K > 1/\delta_1$,

$$c^2 \leq F_2(\delta_1) + \delta_1 \log(\#\Sigma) + \delta_1(\log 3 + 1).$$

We conclude that there is an $F_3$ with

$$c \leq F_3(\delta_1).$$

Let $S_t = S_0 \backslash S_1$ and we have:

1) $p_t(S_t) \geq 1 - F_2(\delta_1) - F_3(\delta_1)$,

2) $\#S_t \leq \#S_0 \leq 2^{(h+\delta_1)\#K}$, and

3) for $\eta \in S_t$, $p_t(\eta) \leq 2^{-(h-F_3(\delta_1))\#K}$.

That is to say for $F_4(z) = F_2(z) + F_3(z) + z$, $p_t$ is $(h, F_4(\delta_1))$-flat. For $\delta > 0$ choose $\delta_1$ so that $F_4(\delta_1) \leq \delta$ to give the result. □

Of course the disintegration of the distribution envisioned in this result will arise from conditioning on a sub-$\sigma$-algebra.

We now start to work with conditional entropy. Let $\mathcal{A}$ be a sub $\sigma$-algebra and $P$ a finite partition. For our purposes it is best to think of this partition as a map $P : X \to \Sigma$ where $\Sigma$ is a finite labeling set. The partition in the usual sense is then the collection of sets $\{P^{-1}(s)\}_{s \in \Sigma}$. As described earlier by $E(P|\mathcal{A})$ we mean the probability vector-valued function $(E(\mathbf{1}_{P^{-1}(s)}|\mathcal{A}))_{s \in \Sigma}$. By $L(P|\mathcal{A})$ we mean the conditional entropy as an $\mathcal{A}$-measurable function

$$L(P|\mathcal{A})(x) = H(E(P|\mathcal{A})(x))$$

and now the conditional entropy of $P$ given the $\sigma$-algebra $\mathcal{A}$ is defined as

$$h(P|\mathcal{A}) = \int_X L(P|\mathcal{A}) \, d\mu.$$

LEMMA 4.4. *For $T$ an ergodic action of $\Gamma$, $P$ and $Q$ two finite partitions, and any $\delta > 0$ if $F$ is sufficiently invariant, then for all but $\delta$ in measure of $X$*

$$E(\bigvee_{\gamma \in F} T_{\gamma^{-1}}(P) | \bigvee_{\gamma \in F} T_{\gamma^{-1}}(Q))$$

*is $(\delta, h)$-flat where $h = h(T, P \vee Q) - h(T, Q)$.*

*Proof.* From Corollary 3.0.26 of [4] which is a simple extension of the Shannon-McMillan theorem of [7] (§II.4 #5), if $F$ is sufficiently invariant, then for all but $(\delta/2)^2$ in measure of the atoms $\eta$ of $\bigvee_{\gamma \in F} T_{\gamma^{-1}}(P \vee Q)$ we have

1) $$\mu(\eta) = 2^{-(h(T, P \vee Q) \pm (\delta/2)^2) \#F}$$

and for all but $\delta/2$ in measure of the atoms $\eta' \in \bigvee_{\gamma \in F} T_{\gamma^{-1}}(Q)$,

2) $$\mu(\eta') = 2^{-(h(T, Q) \pm \delta/2) \#F}.$$

Consider those atoms $\eta'$ satisfying 2) and all but a fraction $\delta/2$ in measure covered by atoms $\eta$ satisfying 1). For $x$ in such an $\eta'$ and also in such an $\eta$ we compute

$$E(\eta | \bigvee_{\gamma \in F} T_{\gamma^{-1}}(Q))(x) = \frac{\mu(\eta)}{\mu(\eta')} = 2^{-(h \pm \delta) \#F};$$

hence

$$E(\bigvee_{\gamma \in F} T_{\gamma^{-1}}(P) | \bigvee_{\gamma \in F} T_{\gamma^{-1}}(Q))$$

is $(\delta, h)$-flat. □

Our next problem is to show that if $F$ is sufficiently invariant then

$$h(\bigvee_{\gamma \in F} T_{\gamma^{-1}}(P) | \bigvee_{\gamma \in \Gamma} T_{\gamma^{-1}}(Q))/\#F$$

cannot be much smaller than

$$h(\bigvee_{\gamma \in F} T_{\gamma^{-1}}(P) | \bigvee_{\gamma \in F} T_{\gamma^{-1}}(Q))/\#F.$$

THEOREM 4.5. *For $T$ an ergodic action of $\Gamma$, $P$ and $Q$ two finite partitions and $\delta > 0$, if $F$ is sufficiently invariant then*

$$\frac{1}{\#F} h(\bigvee_{\gamma \in F} T_{\gamma^{-1}}(P) | \bigvee_{\gamma \in \Gamma} T_{\gamma^{-1}}(Q)) \geq h(T, P \vee Q) - h(T, Q) - \delta.$$

*Proof.* Suppose $F \subseteq F'$ and notice that the calculation

$$\frac{1}{\#F} h(\bigvee_{\gamma \in F} T_{\gamma^{-1}}(P) | \bigvee_{\gamma \in F'} T_{\gamma^{-1}}(Q))$$

can only decrease as $F'$ increases as a set. Suppose the result is false, that is to say fails to hold for some $\delta_0 > 0$. This would mean we will be able to find arbitrarily invariant $F$'s for which it fails. In particular for any $\delta_1 > 0$ we will be able to find $F_1, \ldots F_N$ which $\delta_1$-quasi tile all sufficiently invariant $F$'s and a single $F'$ finite with

$$\frac{1}{\#F_i} h(\bigvee_{\gamma \in F_i} T_{\gamma^{-1}}(P) | \bigvee_{\gamma \in F'} T_{\gamma^{-1}}(Q)) < h(T, P \vee Q) - h(T, Q) - \delta_0/2.$$

Choose $F$ to also be $[F', \delta_1]$-invariant. Now we take a $\delta_1$-quasi tiling of $F$ by $F_1, \ldots, F_N$ with centers $c_{i,j}$, $i = 1, \ldots, N$ and $j = 1, \ldots, k(i)$. Delete from the list of centers all $c_{i,j}$ with $F_i c_{i,j}$ not contained in $F$. The remainder of the tiles will still $2\delta_1$-quasi tile $F$.

We calculate with this an upper bound for

$$\frac{1}{\#F} h(\bigvee_{\gamma \in F} T_{\gamma^{-1}}(P) | \bigvee_{\gamma \in F} T_{\gamma^{-1}}(Q))$$

by successively adding on a conditional entropy of a span over a new tile of indices $F_i c_{i,j}$ until they are exhausted and then adding a maximal bound for the added entropy of the remaining indices. This gives the inequalities

$$\frac{1}{\#F} h(\bigvee_{\gamma \in F} T_{\gamma^{-1}}(P) | \bigvee_{\gamma \in F} T_{\gamma^{-1}}(Q))$$
$$\leq \frac{1}{\#F} \left( \sum_{i,j} h(\bigvee_{\gamma \in F_i} T_{\gamma^{-1}}(P) | \bigvee_{\gamma \in F'} T_{\gamma^{-1}}(Q)) \#F_i \right) + 2\delta_1 \log(\#\Sigma_P)$$
$$\leq \frac{\sum_{i,j} F_i}{\#F} (h(T, P \vee Q) - h(T, Q) - \delta_0/2) + 2\delta_1 \log(\#\Sigma_P)$$
$$\leq (1 + \delta_1)(h(T, P \vee Q) - h(T, Q)) + 2\delta_1 \log(\Sigma_P) - \delta_0/2.$$

If $F$ is sufficiently invariant we conclude from Lemma 4.4 that

$$(1) \quad h(T, P \vee Q) - h(T, P) - \delta_1 \leq (1 + \delta_1)(h(T, P \vee Q) - h(T, Q))$$
$$+ 2\delta_1 \log(\Sigma_P) - \delta_0/2$$

for all values $\delta_1 > 0$. But this conflicts with $\delta_0 > 0$. □

COROLLARY 4.6. *For any $T$-invariant sub-$\sigma$-algebra $\mathcal{A}$ and $\delta > 0$, if $F$ is sufficiently invariant, for all but $\delta$ in measure of the $x \in X$, the conditional measure*

$$E(\bigvee_{\gamma \in F} T_{\gamma^{-1}}(P) | \mathcal{A})(x)$$

*is $(h(T, P|\mathcal{A}), \delta)$-flat.*

*Proof.* Fix a partition $Q$. As the invariance of $F$ improves the measures
$$E(\bigvee_{\gamma \in F} T_{\gamma^{-1}}(P) | \bigvee_{\gamma \in F} T_{\gamma^{-1}}(Q))(x)$$
become ever flatter on ever more of $X$, and the amount the entropy can decrease when further conditioning from $\bigvee_{\gamma \in F} T_{\gamma^{-1}}(Q)$ to $\bigvee_{\gamma \in \Gamma} T_{\gamma^{-1}}(Q)$ becomes ever smaller. Thus Theorem 4.3 gives the result for $\mathcal{A}$ of the form $\bigvee_{\gamma \in \Gamma} T_{\gamma^{-1}}(Q)$. Any invariant sub-$\sigma$-algebra is an increasing span of such algebras $\mathcal{A}_i$. The conditional entropies $h(\bigvee_{\gamma \in F} T_{\gamma^{-1}}(P) | \mathcal{A}_i)/\#F$ decrease in $i$ converging to
$$h(\bigvee_{\gamma \in F} T_{\gamma^{-1}}(P) | \mathcal{A})/\#F.$$
This latter, as $F$ becomes larger and more invariant, decreases to a value we call $h(T, P | \mathcal{A})$.

Thus for any $\delta > 0$, once $i$ is large enough and $F$ sufficiently invariant,
$$\frac{1}{\#F} h(\bigvee_{\gamma \in F} T_{\gamma^{-1}}(P) | \mathcal{A}) \geq \frac{1}{\#F} h(\bigvee_{\gamma \in F} T_{\gamma^{-1}}(P) | \mathcal{A}_i) - \delta.$$
Once $F$ is sufficiently invariant for all but $\delta$ in measure of the $x \in X$, we already know $E(\bigvee_{\gamma \in F} T_{\gamma^{-1}}(P) | \mathcal{A}_i)$ will be $(h(T, P \vee Q_i) - h(T, Q_i), \delta)$-flat. Once more Theorem 4.3 gives the result. □

COROLLARY 4.7. *As the invariance of a set $F$ increases to infinity,*
$$\frac{1}{\#F} L(\bigvee_{\gamma \in F} T_{\gamma^{-1}}(P) | \mathcal{A}) \to h(T, P | \mathcal{A}).$$

*Proof.* Notice that
$$\frac{1}{\#F} L(\bigvee_{\gamma \in F} T_{\gamma^{-1}}(P) | \mathcal{A})(x) \leq \log(\#\Sigma_P).$$
The increasing flatness of the measures $E(\bigvee_{\gamma \in F} T_{\gamma^{-1}}(P) | \mathcal{A})$ now gives the result. □

The corollary above can be obtained quite easily as well from Theorem 3.2 of [9] where it is shown that $I(\bigvee_{\gamma \in F} T_{\gamma^{-1}}(P) | \mathcal{A})/\#F$ converges in $L_1$ as the invariance of $F$ increases.

Our next step is to show that conditional entropy can be found by examining it on sufficiently invariant castles. This is the critical step in proving Theorem 2.6. Our approach is parallel to that of our proof of the ergodic theorem on castles, we will show that if it holds on *some* sequence of castles of growing invariance then it must hold on all such. To get the existence of a sequence we use the same trick as before, we add on another coordinate and note that then it is implied by the previous corollary.

To begin, for any finite subset $H \subset \Gamma$ and $A \subset X$ with the pair forming a tower we can calculate a local conditional entropy on the tower $H \times A$ as

$$h_{(H,A)}(T, P | \mathcal{A}) = \int_A L(\bigvee_{\gamma \in H} T_{\gamma^{-1}}(P) | \mathcal{A}).$$

Notice we have not normalized this by either $\mu(A)$ or by $\#H$. The reason is easy to understand from the inequalities

$$h_{(H,A)}(T, P | \mathcal{A}) \leq \int \#H \log n \, d\mu \leq \log n \#H \mu(A)$$

where $P$ is an $n$-set partition. That is to say this tower entropy is automatically bounded by the measure of the tower image times $\log n$. For a castle $\mathcal{T} = \{H_i, A_i\}$ we can define a conditional entropy on the castle now as

$$h_{\mathcal{T}}(T, P | \mathcal{A}) = \sum_i h_{(h_i, A_i)}(T, P | \mathcal{A})$$

and again we see that this tower entropy will be bounded by the measure of the castle image times $\log n$.

LEMMA 4.8. *Suppose $T$ is a free and ergodic action of the countable and amenable group $\Gamma$, $\mathcal{A}$ is an invariant sub-$\sigma$-algebra and $P$ is an $n$-set partition. Suppose $\mathcal{T} = \{H_i, A_i\}$ is a castle covering all but $\delta \geq 0$ of $X$ with $A_i \in \mathcal{A}$. Given any $\varepsilon > 0$ there is an $N$, so that for all castles $\mathcal{T}' = \{K_j, B_j\}$ that are $N$-invariant with $B_j \in \mathcal{A}$ we will have*

$$h_{\mathcal{T}'}(T, P | \mathcal{A}) \leq h_{\mathcal{T}}(T, P | \mathcal{A}) + \delta \log n + \varepsilon.$$

*Proof.* We will abuse language a bit here and use $\mu(\mathcal{T})$ and $\mu(\mathcal{T}')$ to represent measures of the two castle images. We assume now that the partition $P$ and the castle $\mathcal{T}$ are given. For a castle $\mathcal{T}' = \{K_j, B_j\}$ and $x \in B_j$ define

(2) $\qquad d(x) = \{\gamma \in K_j : \text{ there exists an } i \text{ with}$
$$\gamma = \gamma_2 \gamma_1, \, T_{\gamma_1}(x) \in A_i, \, \gamma_2 \in H_i \text{ and } H_i \gamma_1 \subseteq K_j\}$$

that is to say the collection of elements in $K_j$ which lie in a slice through the $\mathcal{T}$ castle which lies completely inside this slice through the $\mathcal{T}'$-castle.

We note that if $\mathcal{T}'$ is sufficiently invariant, by the ergodic theorem on towers and this invariance we will obtain

$$\sum_j \int_{B_j} \left| \#d(x) - \#K_j \mu(\mathcal{T}) \right| d\mu < \frac{\varepsilon}{2 \log n}.$$

Define a function $h$ on $X$ by

$$h(x) = \begin{cases} \frac{1}{\#H_i} L(\bigvee_{\gamma \in H_i} T_{\gamma^{-1}}(P) | \mathcal{A})(x') & \text{if } x = T_\gamma(x'), \, \gamma \in H_i, \, x' \in A_i \\ 0 & \text{otherwise.} \end{cases}$$

Notice that $h(x) \leq \log n$ and
$$\int h(x)\, d\mu = h_{\mathcal{T}}(T, P|\mathcal{A}).$$

Using the ergodic theorem on towers, once $\mathcal{T}'$ is sufficiently invariant we will obtain
$$\sum_j \int_{B_j} \left| \sum_{\gamma \in K_j} h(T_\gamma(x)) - \#K_j h_{\mathcal{T}}(T, P|\mathcal{A}) \right| d\mu < \varepsilon/2$$

and in particular
$$\sum_j \int_{B_j} \sum_{\gamma \in K_j} h(T_\gamma(x))\, d\mu < h_{\mathcal{T}}(T, P|\mathcal{A}) + \varepsilon/2.$$

For $x \in B_j$ for a first estimate
$$L(\bigvee_{\gamma \in K_j} T_{\gamma^{-1}}(P)|\mathcal{A})(x) \leq L(\bigvee_{\gamma \in d(x)} T_{\gamma^{-1}}(P)|\mathcal{A})(x) + (\#K_j - \#d(x))\log n.$$

Let $d_i(x) = \{\gamma \in d(x) : T_\gamma(x) \in A_i\}$ and now
$$d(x) = \bigcup_i \bigcup_{\gamma_1 \in d_i(x)} \bigcup_{\gamma_2 \in H_i} T_{\gamma_2 \gamma_1}(x)$$

a disjoint union.

Choose some ordering for each of the $d_i(x) = \{\gamma_{ik}\}$ and set
$$C_{i,k} = \{\gamma \in d(x) : \gamma = \gamma_1 \gamma_{i'k'} \text{ where } \gamma_1 \in A_{i'},\ i' \leq i \text{ and } k' < k\}.$$

What this has done is simply to give an ordering to the slices through $\mathcal{T}$ that lie within the slice of $\mathcal{T}'$ passing through $x$ with $C_{i,k}$ consisting of those elements in the slice that occur earlier in the order than the one through $\gamma_{i,k}$.

As both $\mathcal{T}$ and $\mathcal{T}'$ are $\mathcal{A}$-measurable we compute that
$$\begin{aligned}
I(\bigvee_{\gamma \in d(x)} T_{\gamma^{-1}}(P)|\mathcal{A})(x) &= \sum_{i,k} L(\bigvee_{\gamma \in A_i} T_{\gamma^{-1}}(P)|\bigvee_{\gamma \in C_{i,k}} T_{\gamma^{-1}}(P) \vee \mathcal{A})(T_{\gamma_{i,k}}(x)) \\
&\leq \sum_{i,k} L(\bigvee_{\gamma \in A_i} T_{\gamma^{-1}}(P)|\mathcal{A})(T_{\gamma_{i,k}}(x)) \\
&= \sum_{\gamma \in d(x)} h(T_\gamma(x)).
\end{aligned}$$

We now complete our work by combining our estimates to give
$$\begin{aligned}
h_{\mathcal{T}'}(T, P|\mathcal{A}) &= \sum_j \int_{B_j} L(\bigvee_{\gamma \in K_j} T_{\gamma^{-1}}(P)|\mathcal{A})\, d\mu \\
&\leq \sum_j \int_{B_j} L(\bigvee_{\gamma \in d(x)} T_{\gamma^{-1}}(P)|\mathcal{A})(x)\, d\mu(x) \\
&\quad + \sum_j \int_{B_j} (\#K_j - \#d(x))\log n\, d\mu
\end{aligned}$$

$$\leq \sum_j \int_{B_j} \sum_{\gamma \in d(x)} h(T_\gamma(x))\, d\mu$$

$$+ \left(\sum_j \#K_j \mu(B_j)(1-\mu(\mathcal{T}))\right) \log n + \varepsilon/2$$

$$\leq \sum_j \int_{B_j} \sum_{\gamma \in K_j} h(T_\gamma(x))\, d\mu + \delta \log n + \varepsilon/2$$

$$\leq h_\mathcal{T}(T, P|\mathcal{A}) + \delta \log n + \varepsilon. \qquad \square$$

COROLLARY 4.9. *For $T$ an ergodic action of the countable and amenable group $\Gamma$, $\mathcal{A}$ a $T$-invariant sub-$\sigma$-algebra and $P$ an $n$-set partition, then for any $\varepsilon > 0$ there is an $N$ and a $\delta$ so that if $\mathcal{T} = \{H_i, A_i\}$ is a castle that is $N$-invariant and covers all but $\delta$ of $X$, then*

$$h(T, P|\mathcal{A}) + \delta \log n + \varepsilon < h_\mathcal{T}(T, P|\mathcal{A}) < h(T, P|\mathcal{A}) + \varepsilon.$$

*Proof.* We split the proof into two pieces. First, suppose that $T$ possesses a sequence of castles $\mathcal{T}_i$ which are $\mathcal{A}$-measurable, whose invariance tends to infinity, and for which $\mu(\mathcal{T}_i) \to 1$ and

$$h_{\mathcal{T}_i}(T, P|\mathcal{A}) \to h(T, P|\mathcal{A}).$$

Choose a $\mathcal{T}_i$ covering at least $\varepsilon/(3 \log n)$ of $X$ and so that

$$h_{\mathcal{T}_i}(T, P|\mathcal{A}) < h(T, P|\mathcal{A}) + \varepsilon/3$$

and using $\varepsilon/3$ in the previous lemma we conclude that once $\mathcal{T}$ is sufficiently invariant and $\mathcal{A}$-measurable we will have

$$h_\mathcal{T}(T, P|\mathcal{A}) < h_{\mathcal{T}_i}(T, P|\mathcal{A}) + 2\varepsilon/3 < h(T, P|\mathcal{A}) + \varepsilon.$$

Now fixing such a $\mathcal{T}$ and again using the previous lemma with $\varepsilon/2$, choose $\mathcal{T}_j$ with $h_{\mathcal{T}_j}(T, P|\mathcal{A}) > h(T, P|\mathcal{A}) - \varepsilon/2$ and giving now

$$h(T, P|\mathcal{A}) < h_{\mathcal{T}_j}(T, P|\mathcal{A}) + \varepsilon/2 < h_\mathcal{T}(T, P|\mathcal{A}) + \varepsilon/2 + \delta \log n + \varepsilon/2$$

giving the result.

Where will we find the castles $\mathcal{T}_i$? We use the same trick as we did for the ergodic theorem on towers. Construct a direct product $T \times T_2$ that is still ergodic. We leave $P$ fixed but replace $\mathcal{A}$ with $\mathcal{A}_1 = \mathcal{A} \vee \mathcal{B}_2$. For any tower $\{H, A\}$ that is $\mathcal{B}_2$-measurable,

$$h_{(H,A)}(T \times T_2 | \mathcal{A}_1) = h(\bigvee_{\gamma \in H} T_{\gamma^{-1}}(P)|\mathcal{A})/\#H.$$

From Corollary 4.7 we know that as the invariance of $H$ grows,

$$h(\bigvee_{\gamma \in H} T_{\gamma^{-1}}(P)|\mathcal{A})/\#H$$

will converge to $h(T, P|\mathcal{A})$. Thus if we choose a sequence of castles $\mathcal{T}_i$ which are $\mathcal{B}_2$-measurable, whose invariance grows and which cover more and more of $X \times X_2$ we will have the sequence of castles we want. We now obtain the conclusion for all castles in the product action that are $\mathcal{A}_1$-measurable and in particular for those just in the first coordinate and $\mathcal{A}$-measurable. □

*Proof of Theorem* 2.6. From Theorem 2.11 we know that for any sequence of castles $\mathcal{T}_i$ for the action of $T$ whose invariance grows to infinity, the image castles $\mathcal{T}_i'$ for the action of $U$ will have the same property. If the $\mathcal{T}_i$ are $\mathcal{A}$-measurable, as the orbit change is also $\mathcal{A}$-measurable, the castles $\mathcal{T}_i'$ will again be $\mathcal{A}$-measurable. Hence we must have both

$$\lim_{i \to \infty} h_{\mathcal{T}_i}(T, P|\mathcal{A}) = h(T, P|\mathcal{A})$$

and

$$\lim_{i \to \infty} h_{\mathcal{T}_i'}(U, P|\mathcal{A}) = h(U, P|\mathcal{A}).$$

As the orbit change is $\mathcal{A}$-measurable, for almost every $x \in X$ and finite set $H \subseteq \Gamma$, letting $H' = \alpha(x, H)$,

$$L(\vee_{\gamma \in H} T_{\gamma^{-1}}(P)|\mathcal{A})(x) = L(\vee_{\gamma' \in H'} U_{\gamma'^{-1}}(P)|\mathcal{A})(x);$$

hence for all $\mathcal{T}_i$ we have

$$h_{\mathcal{T}_i}(T, P|\mathcal{A}) = h_{\mathcal{T}_i'}(U, P|\mathcal{A})$$

which completes the result. □

As promised in the introduction, we will finish by showing that in certain special cases it is easy to see why the relative Pinsker algebra of a direct product over its first coordinate is simply the span of the first coordinate and the Pinsker algebra of the second. In particular our argument will hold when the first coordinate action is Bernoulli.

THEOREM 4.10. *Suppose $T_1$ and $T_2$ are two ergodic and measure-preserving actions of $\Gamma$ on the spaces $(X_1, \mathcal{B}_1, \mu_1)$ and $(X_2, \mathcal{B}_2, \mu_2)$ respectively. Let $\Pi(T_i)$ represent the Pinsker algebra of $T_i$ and $\mathcal{E}$ represent the relative Pinsker algebra of the direct product $T_1 \times T_2$ with respect to the first coordinate algebra $\mathcal{B}_1$. Suppose that $C(T_1)$, the group of all measure-preserving bijections of $X_1$ commuting with all elements of $T_1$, acts ergodically on $X_1$. We conclude that*

$$\mathcal{E} = \mathcal{B}_1 \times \Pi(T_2).$$

*Proof.* Denote by $\mathcal{E}_{x_1}$ the trace of $\mathcal{E}$ on the fiber $x_1 \times X_2$. Up to $\mu_2$-null sets we can identify $\mathcal{E}_{x_1}$ with the conditional expectation operator from $L_2(X_2, \mathcal{B}_2, \mu_2)$ to $L_2(X_2, \mathcal{E}_{x_1}, \mu_2)$ which we write as $E_{x_1}$. Using the strong operator topology on these conditional expectations, we get a map from $X_1$ into a subset of the closed subset of conditional expectation operators in this Polish space. Notice that for any function $f \in L_2(X_2, \mathcal{B}_2, \mu_2)$ we can lift $f$ to $X_1 \times X_2$ as $f(x_1, x_2) = f(x_2)$ and we have

$$E_{x_1}(f)(x_2) = E(f|\mathcal{E})(x_1, x_2).$$

This shows that for any $f$ the evaluation of $E_{x_1}(f)(x_2)$ is jointly measurable in both variables and hence the map $E_{x_1}$ is a measurable map to this Polish (in particular separable metric) space. If $S$ is an element of $C(T_1)$ the fact that $S \times \text{id}(\mathcal{E}) = \mathcal{E}$ tells us that $E_{Sx_1} = E_{x_1}$ for $\mu_1$-a.e. $x_1$. The ergodicity of $C(T_1)$ now implies that the operators $E_{x_1}$ must be $\mu_1$-almost surely constant, and hence the algebras $\mathcal{E}_{x_1}$ must be almost surely a constant $\mathcal{E}_0$. Since $\mathcal{B}_1$ is contained in $\mathcal{E}$ we conclude that $\mathcal{E} = \mathcal{B}_1 \times \mathcal{E}_0$ and it is an observation that $\mathcal{E}_0 = \Pi(T_2)$ as any set in $\mathcal{B}_2$ with zero conditional entropy over $\mathcal{B}_1$ must have zero entropy for the action $T_2$. □

For our purposes here all we need is a single example of such a map $T_1$ for which we can guarantee the $T_1 \times T_2$ will remain ergodic for all ergodic $T_2$. Here is perhaps the simplest such example. Let $X = \{0,1\}^\Gamma$ be the full 2-shift over $\Gamma$ with the *left* action, i.e.

$$T_\gamma(x_1)(\gamma') = x_1(\gamma^{-1}\gamma')$$

and let $\mu_1$ be independent and identically distributed measure $(1/2, 1/2)^\Gamma$. Notice that $C(T_1)$ will contain the *right* action of $\Gamma$

$$S_\gamma(x_1)(\gamma') = x_1(\gamma'\gamma_1).$$

The action $S$ is just as ergodic as that of $T$.

COROLLARY 4.11. *If $T_1$ acting on $(X_1, \mathcal{B}_1, \mu_1)$ is Bernoulli and $T_2$ acting on $(X_2, \mathcal{B}_2, \mu_2)$ is cpe then the product action $T_1 \times T_2$ is relatively cpe over $X_1$.*

## 5. Final comments

Recent work in this area has not focused on proving that completely positive entropy implies multiple mixing, but rather on the characterizing the unitary theory of such systems. In particular the goal is to show that the spectral type of any system of a completely positive entropy action of $\Gamma$ is that of the shift on $\ell^2$ of the group $\Gamma$. This has been shown for $\mathbb{Z}^n$ by B. Kaminsky and for the discrete rationals by J.-P. Thouvenot. Is it possible to approach this question using the orbit equivalence notions presented here? That is to say can one pull the relativized spectral structure through an orbit equivalence to show that the result is simply the transference of the classical fact for $\mathbb{Z}$?

Having restricted our attention here to countable amenable groups one must ask whether the methods here will apply to continuous such groups. The answer is certainly going to be yes. For groups that possess a cocompact countable subgroup the result on the subgroup together with some coding approximations will give the result for the larger group. For groups which do not possess such a subgroup the path is a bit trickier. Remember that central to our proof was the taking of a direct product with a standard second action. We once more do this but in this standard action we also select a section. The induced equivalence relation on this section will be hyperfinite of type $II_1$, that is to say orbit-equivalent to a measure-preserving action of $\mathbb{Z}$. It is through this orbit equivalence that we will transfer the relative cpe property and mixing properties to obtain our result. This work will be presented separately.

We end with a tantalizing open question. An area of the ergodic theory of discrete amenable groups that remains very cloudy is the nature of Følner sequences along which the pointwise ergodic theory holds for all $f \in L^1$. We pose the following question: Suppose $T$ is some fixed action of the countable and discrete amenable group $\Gamma$ and $U$ is an action of $\mathbb{Z}$ with the same orbits. For a.e. $x \in X$ Lemma 3.10 tells us the sequence of sets $F_n(x) = \{\gamma \in \Gamma : T_\gamma(x) = U_i(x), -n \leq i \leq n\}$ is a Følner sequence in $\Gamma$. Is it the case that for a.e. $x \in X$ this sequence is universally good for the pointwise ergodic theorem? What is easily seen from our transference is that for any action $S$ of $\Gamma$ on a space $Y$ and any $f \in L^1(Y)$ that for almost every $x \in X$, the sequence $F_n(x)$ is good for the pointwise ergodic theorem for the $S$ action on $f$.

University of Maryland, College Park, MD
*E-mail address*: djr@math.umd.edu

Mathematics Institute, Hebrew University of Jerusalem, Jerusalem, Israel
*E-mail address*: weiss@math.huji.ac.il